 \documentclass[final,authoryear,1p,times]{elsarticle}
\usepackage{amssymb}
\usepackage{amsthm}

\usepackage{graphicx}
\usepackage{subfigure}

\newtheorem{theorem}{Theorem}[section]

\newtheorem{conjecture}[theorem]{Conjecture}

\journal{}

\begin{document}

\begin{frontmatter}



\title{Conversion of matrix weighted rational B\'{e}zier curves to rational B\'{e}zier curves}


\author{Xunnian Yang\corref{cor1}}
\date{}
%\ead[cor1]{Tel.: +86 571 87951609; fax: +86 571 87951428.}
\ead{yxn@zju.edu.cn}

\address{School of mathematical sciences, Zhejiang University, Hangzhou 310027, China}

\begin{abstract}
Matrix weighted rational B\'{e}zier curves can represent complex curve shapes using small numbers of control points and clear geometric definitions of matrix weights. Explicit formulae are derived to convert matrix weighted rational B\'{e}zier curves in 2D or 3D space to rational B\'{e}zier curves. A method for computing the convex hulls of matrix weighted rational B\'{e}zier curves is given as a conjecture.
\end{abstract}


\end{frontmatter}

% \linenumbers


%%%%%%%%%%%%%%%%%%%%%%%%%%%%%%%%%%%%%%%%%%%%%%%%%%%%%%%%%%%%%%%%%%%%%%%%%%%%%%%%%%%%%%%%%%%
%%% Sectioin 1                                                                          %%%
%%%%%%%%%%%%%%%%%%%%%%%%%%%%%%%%%%%%%%%%%%%%%%%%%%%%%%%%%%%%%%%%%%%%%%%%%%%%%%%%%%%%%%%%%%%

\section{Introduction}
\label{Sec:intro}

Rational curves and surfaces are powerful tools for shape design in Computer Aided Geometric Design~\citep{Farin2001CAGDBook}. Besides the real numbers, the weights can be chosen as complex numbers for planar rational curves~\citep{Sanchez-Reyes09ComplexBezier} or matrices for rational curves and surfaces in arbitrary dimension space~\citep{Yang16MatrixRational}. Particularly, matrix weighted rational B\'{e}zier curves can be used to design complex shapes with small numbers of control points while matrix weighted NURBS curves can be used to fit and fair Hermite-type data by proper geometric definitions of the matrix weights~\citep{Yang18MatrixNURBS}.

Suppose that $P_0, P_1, \ldots, P_n$ are a sequence of points lying in $\mathbb{R}^d$ and $M_i\in\mathbb{R}^{d\times d}$, $i=0,1,\ldots,n$, are a set of weight matrices. A matrix weighted rational B\'{e}zier curve of degree $n$ is defined by
\begin{equation}\label{Eqn:Q(t)_Bezier}
Q(t)=\left[\sum_{i=0}^n M_i B_{i,n}(t)\right]^{-1}\sum_{i=0}^n M_iP_i B_{i,n}(t), \ \ \ t\in[0,1],
\end{equation}
where $B_{i,n}(t)={{n!}\over{i!(n-i)!}}t^i(1-t)^{n-i}$ are the Bernstein basis functions.

For ease of shape control, the weight matrices for a matrix weighted rational B\'{e}zier curve can be defined by normal vectors~\citep{Yang16MatrixRational} or tangent vectors specified at the control points~\citep{Yang18MatrixNURBS}. Suppose that $\mathbf{n}_i$, $i=0,1,\ldots,n$, are unit normal vectors and $\omega_i>0$, $\mu_i>-1$, $i=0,1,\ldots,n$, are real numbers. The weight matrices for matrix weighted rational B\'{e}zier curves with point-normal control pairs are given by
\begin{equation}\label{Eqn:M_i by point_normal}
M_i=\omega_i(I+\mu_i\mathbf{n}_i\mathbf{n}_i^T), \ \ i=0,1,\ldots,n,
\end{equation}
where $I$ is the identity matrix of order $d$ and $T$ means the transpose of a column vector.
If a set of unit tangents $\mathbf{t}_i$ have been specified at the control points, the weight matrices are then given by
\begin{equation}\label{Eqn:M_i by point_tangent}
M_i=\omega_i[I+\mu_i(I-\mathbf{t}_i\mathbf{t}_i^T)], \ \ i=0,1,\ldots,n.
\end{equation}
If the weight matrices of a matrix weighted rational B\'{e}zier curve are defined by Equation~(\ref{Eqn:M_i by point_tangent}), the shape of the curve will be controlled efficiently by control points and tangent lines passing through the points.

In~\cite{Yang16MatrixRational} we have proven that matrix weighted rational B\'{e}zier curves are actually the conventional rational B\'{e}zier curves. In the following two sections we will derive explicit formulae for converting matrix weighted rational B\'{e}zier curves in 2D or 3D space into standard rational B\'{e}zier curves. Finally, a conjecture for computing the convex hulls of matrix weighted rational B\'{e}zier curves by converting them to rational B\'{e}zier curves will be given.

%Though matrix weighted rational B\'{e}zier curves have nice properties, they do not lie in the convex hulls of their control points generally. The convex hulls of matrix weighted rational B\'{e}zier curves can be computed. The conversion can also help to transfer the curves in different modeling systems.

%%%%%%%%%%%%%%%%%%%%%%%%%%%%%%%%%%%%%%%%%%%%%%%%%%%%%%%%%%%%%%%%%%%%%%%%%%%%%%%%%%%%%%%%%%%
%%% Sectioin 2                                                                          %%%
%%%%%%%%%%%%%%%%%%%%%%%%%%%%%%%%%%%%%%%%%%%%%%%%%%%%%%%%%%%%%%%%%%%%%%%%%%%%%%%%%%%%%%%%%%%

\section{Convert a matrix weighted rational B\'{e}zier curve to a rational B\'{e}zier curve in 2D}
\label{Sec:ConversionIn2D}

Let $M(t)=\sum_{i=0}^n M_i B_{i,n}(t)$ and $M^*(t)$ be the adjoint matrix of $M(t)$. Then the matrix weighted rational B\'{e}zier curve given by Equation~(\ref{Eqn:Q(t)_Bezier}) can be reformulated as
\begin{equation}\label{Eqn:Q(t)_rational Bezier}
Q(t)=\frac{M^*(t)\sum_{i=0}^n M_iP_i B_{i,n}(t)}{\det(M(t))}, \ \ \ t\in[0,1].
\end{equation}
To convert the matrix weighted rational B\'{e}zier curve to a standard rational B\'{e}zier curve, we should then represent both the numerator and the denominator of Equation~(\ref{Eqn:Q(t)_rational Bezier}) by Bernstein basis functions.

Assume that $Q(t)$ is a matrix weighted rational B\'{e}zier curve in 2D and the matrix function $M(t)$ is given by
\[
\begin{array}{lll}
M(t)&=&\left(
    \begin{array}{cc}
        a(t) & b(t) \\
        c(t) & d(t)
    \end{array}
    \right) \\
    &=& \left(
    \begin{array}{cc}
        \sum_{i=0}^n a_i B_{i,n}(t) & \sum_{i=0}^n b_i B_{i,n}(t) \\
        \sum_{i=0}^n c_i B_{i,n}(t) & \sum_{i=0}^n d_i B_{i,n}(t)
    \end{array}
    \right).
   % &=& \sum_{i=0}^n\left(
   % \begin{array}{cc}
   %     a_i & b_i \\
   %     c_i & d_i
   % \end{array}
   % \right)B_{i,n}(t).
\end{array}
\]
By the product of polynomials in Bernstein form~\citep{Farouki&Rajan88Algorithmsforpolynomials}, the determinant of the matrix function is obtained as
\begin{equation}\label{Eqn:det(M_2(t))}
\begin{array}{lll}
\det(M(t))&=&a(t)d(t)-b(t)c(t)\\
          &=&\sum_{k=0}^{2n}\omega_k B_{k,2n}(t),
\end{array}
\end{equation}
where $\omega_k=\sum_{i+j=k,0\leq i,j\leq n}\frac{(a_id_j-b_ic_j)C_n^iC_n^j}{C_{2n}^k}$.
The adjoint matrix is obtained as $M^*(t)=\sum_{i=0}^n M_i^* B_{i,n}(t)$, where
\[M_i^*=\left(
    \begin{array}{cc}
        d_i & -b_i \\
        -c_i & a_i
    \end{array}
    \right).
\]

Let $\bar{P}_i=M_iP_i$, $i=0,1,\ldots,n$. The numerator in Equation~(\ref{Eqn:Q(t)_rational Bezier}) can be computed as
\[
\sum_{i=0}^n M_i^* B_{i,n}(t)\sum_{j=0}^n \bar{P}_j B_{j,n}(t)
=\sum_{k=0}^{2n}\omega_k Q_k B_{k,2n}(t)
\]
where
\[
Q_k={\frac{1}{\omega_k}}\sum_{i+j=k,0\leq i,j\leq n} \frac{C_n^i C_n^j M_i^*\bar{P}_j }{C_{2n}^{k}}.
\]
When the weights and control points are computed, the converted rational B\'{e}zier curve is obtained as
\begin{equation}\label{Eqn:Q(t)_2D}
Q(t)=\frac{\sum_{k=0}^{2n}\omega_k Q_k B_{k,2n}(t)}{\sum_{k=0}^{2n}\omega_k B_{k,2n}(t)}.
\end{equation}

\begin{figure}[htb]
  \centering
  \vskip 2cm
  \includegraphics[width=7.2cm]{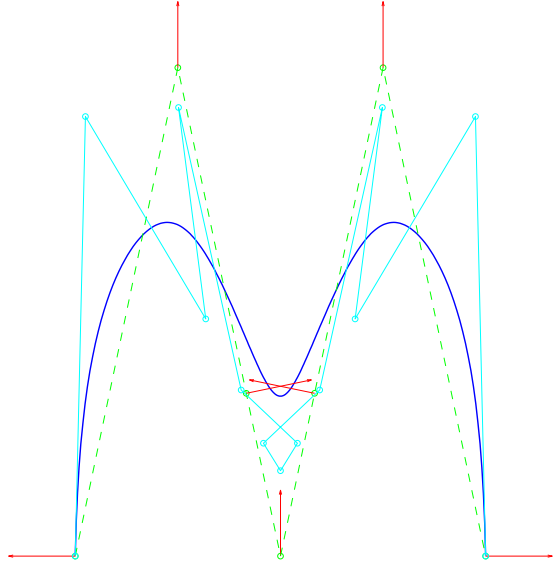}
  \caption{Conversion of a planar matrix weighted rational B\'{e}zier curve with point-normal control pairs into a rational B\'{e}zier curve. }
  \label{Fig:matrix_RaBezier_M}
\end{figure}

In Figure~\ref{Fig:matrix_RaBezier_M}, we design an ``m" like shape by a matrix weighted rational B\'{e}zier curve. By specifying 7 pairs of points and normals consequently, a matrix weighted rational B\'{e}zier curve of degree 6 with point-normal control pairs is defined. The weight matrices are computed by Equation~(\ref{Eqn:M_i by point_normal}) with all $\omega_i=1$ and all $\mu_i=2$ except for $\mu_2=\mu_4=10$. By employing the conversion technique discussed in Section~\ref{Sec:ConversionIn2D}, a rational B\'{e}zier curve of degree 12 is obtained. The control polygon of the converted curve is shown as the solid cyan polyline in Figure~\ref{Fig:matrix_RaBezier_M}. This example demonstrates that small number of control points together with clear geometric definition of matrix weights can help to design a rational curve effectively.

%%%%%%%%%%%%%%%%%%%%%%%%%%%%%%%%%%%%%%%%%%%%%%%%%%%%%%%%%%%%%%%%%%%%%%%%%%%%%%%%%%%%%%%%%%%
%%% Sectioin 3                                                                          %%%
%%%%%%%%%%%%%%%%%%%%%%%%%%%%%%%%%%%%%%%%%%%%%%%%%%%%%%%%%%%%%%%%%%%%%%%%%%%%%%%%%%%%%%%%%%%

\section{Convert a matrix weighted rational B\'{e}zier curve to a rational B\'{e}zier curve in 3D}
\label{Sec:ConversionIn3D}

%This section presents detailed algorithm steps for converting matrix weighted rational B\'{e}zier curves in 3D space to rational B\'{e}zier curves of degree $3n$.

Suppose that $Q(t)$ given by Equation~(\ref{Eqn:Q(t)_Bezier}) is a matrix weighted rational B\'{e}zier curve lying in 3D space and the matrix function $M(t)$ is represented by its elements as
\[
M(t)=\sum_{i=0}^nM_iB_{i,n}(t)=\left(\begin{array}{ccc}
                             a_{11}(t) & a_{12}(t) & a_{13}(t) \\
                             a_{21}(t) & a_{22}(t) & a_{23}(t) \\
                             a_{31}(t) & a_{32}(t) & a_{33}(t)
                           \end{array}
                        \right).
\]
Assume $A_{ij}(t)$ are the corresponding algebraic cofactors of elements $a_{ij}(t)$ of the matrix. We have
\[
\begin{array}{ccc}
A_{11}(t)&=&a_{22}(t)a_{33}(t)-a_{32}(t)a_{23}(t), \\
A_{12}(t)&=&a_{31}(t)a_{23}(t)-a_{21}(t)a_{33}(t), \\
A_{13}(t)&=&a_{21}(t)a_{32}(t)-a_{31}(t)a_{22}(t), \\
A_{21}(t)&=&a_{32}(t)a_{13}(t)-a_{12}(t)a_{33}(t), \\
A_{22}(t)&=&a_{11}(t)a_{33}(t)-a_{31}(t)a_{13}(t), \\
A_{23}(t)&=&a_{31}(t)a_{12}(t)-a_{11}(t)a_{32}(t), \\
A_{31}(t)&=&a_{12}(t)a_{23}(t)-a_{22}(t)a_{13}(t), \\
A_{32}(t)&=&a_{21}(t)a_{13}(t)-a_{11}(t)a_{23}(t), \\
A_{33}(t)&=&a_{11}(t)a_{22}(t)-a_{21}(t)a_{12}(t).
\end{array}
\]
By Equation~(\ref{Eqn:det(M_2(t))}), each $A_{ij}(t)$ can be formulated as a real function of degree $2n$ in terms of Bernstein bases. Then, the adjoint matrix is obtained as
\[
\begin{array}{lll}
M^*(t)&=&\left(\begin{array}{ccc}
                             A_{11}(t) & A_{21}(t) & A_{31}(t) \\
                             A_{12}(t) & A_{22}(t) & A_{32}(t) \\
                             A_{13}(t) & A_{23}(t) & A_{33}(t)
                           \end{array}
                        \right) \\
      &=&\sum_{i=0}^{2n}M_i^* B_{i,2n}(t)
\end{array}
\]
where $M_i^*$, $i=0,1,\ldots,2n$, are the coefficient matrices of order 3.

The real weights of the converted rational B\'{e}zier curve are derived by computing the determinant of the matrix $M(t)$. Assume that $a_{1l}(t)=\sum_{j=0}^n p_{lj}B_{j,n}(t)$, $A_{1l}(t)=\sum_{j=0}^{2n} q_{lj}B_{j,2n}(t)$, $l=1,2,3$. The determinant of the matrix $M(t)$ is formulated as
\begin{equation}\label{Eqn:det_M(t)_3D}
\begin{array}{lll}
\det(M(t))&=& \sum_{l=1}^3a_{1l}(t)A_{1l}(t)  \\
          &=& \sum_{k=0}^{3n} \omega_k B_{k,3n}(t),
\end{array}
\end{equation}
where
\[
\omega_k=\sum_{i+j=k,0\leq i\leq n,0\leq j\leq 2n} \frac{C_n^i C_{2n}^j \sum_{l=1}^3p_{li}q_{lj}}{C_{3n}^{k}}.
\]

Similar to the conversion of matrix weighted rational B\'{e}zier curves in 2D, the control points of the converted rational B\'{e}zier curves in 3D are derived by computing the numerator in Equation~(\ref{Eqn:Q(t)_rational Bezier}). Denote $\bar{P}_i=M_iP_i$, $i=0,1,\ldots,n$. The numerator of the rational B\'{e}zier curve is formulated as
\[
\sum_{i=0}^{2n} M_i^* B_{i,2n}(t)\sum_{j=0}^n \bar{P}_j B_{j,n}(t)
=\sum_{k=0}^{3n}\omega_k Q_k B_{k,3n}(t)
\]
where
\[
Q_k={\frac{1}{\omega_k}}\sum_{i+j=k,0\leq i\leq 2n,0\leq j\leq n} \frac{C_{2n}^i C_n^j M_i^*\bar{P}_j }{C_{3n}^{k}}.
\]
Finally, the converted rational B\'{e}zier curve is
\begin{equation}\label{Eqn:Q(t)_3D}
Q(t)=\frac{\sum_{k=0}^{3n}\omega_k Q_k B_{k,3n}(t)}{\sum_{k=0}^{3n}\omega_k B_{k,3n}(t)}.
\end{equation}

\begin{figure}[htb]
  \centering
  \vskip 2cm
  \includegraphics[width=7.5cm]{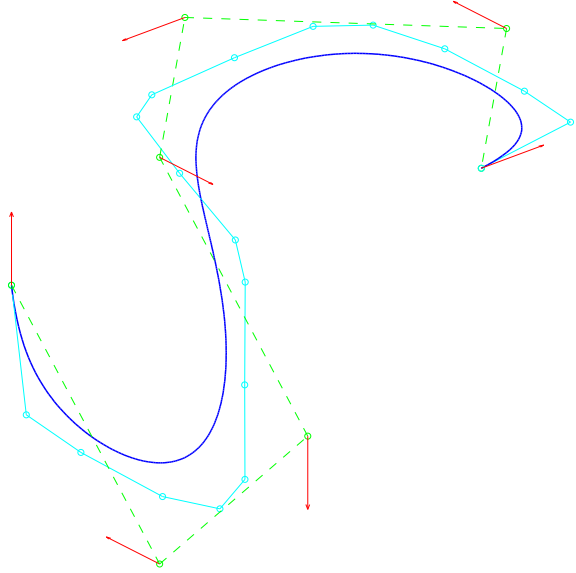}
  \caption{Conversion of a matrix weighted rational B\'{e}zier curve with point-tangent control pairs in 3D space into a rational B\'{e}zier curve. }
  \label{Fig:matrix_RaBezier_S}
\end{figure}

Figure~\ref{Fig:matrix_RaBezier_S} illustrates a spatial matrix weighted rational B\'{e}zier curve of degree 6 with point-tangent control pairs. In particular, the first three point-tangent control pairs and the last three point-tangent control pairs are lying on two planes that are perpendicular with each other while the middle point-tangent control pair lies on the intersection line between the two planes. The weight matrices of the curve are computed by Equation~(\ref{Eqn:M_i by point_tangent}) with all $\omega_i=1$ together with $\mu_0=\mu_6=1.0$, $\mu_1=\mu_5=2.0$ and $\mu_2=\mu_3=\mu_4=4.0$. By applying the techniques discussed in Section~\ref{Sec:ConversionIn3D}, the original curve has been converted to a rational B\'{e}zier curve of degree 18. The solid cyan polygon in Figure~\ref{Fig:matrix_RaBezier_S} illustrates the control polygon of the rational B\'{e}zier curve.

%%%%%%%%%%%%%%%%%%%%%%%%%%%%%%%%%%%%%%%%%%%%%%%%%%%%%%%%%%%%%%%%%%%%%%%%%%%%%%%%%%%%%%%%%%%
%%% Sectioin 4                                                                          %%%
%%%%%%%%%%%%%%%%%%%%%%%%%%%%%%%%%%%%%%%%%%%%%%%%%%%%%%%%%%%%%%%%%%%%%%%%%%%%%%%%%%%%%%%%%%%

\section{Conjecture and closure}
\label{Sec:Conclu}
One main motivation for converting a matrix weighted rational B\'{e}zier curve to a conventional rational B\'{e}zier curve is to compute the convex hull of the curve. A sufficient condition for computing the convex hull of a rational B\'{e}zier curve by its control points is that all the weights are positive. All the examples we have experimented show that the converted rational B\'{e}zier curves do have positive weights, but the result is still difficult to prove at present. We therefore give the result as a conjecture.

\begin{conjecture} Assume a matrix weighted rational B\'{e}zier curve is defined by Equation (\ref{Eqn:Q(t)_Bezier}) with weight matrices given by Equation (\ref{Eqn:M_i by point_normal}) or (\ref{Eqn:M_i by point_tangent}). The real weights computed by Equation (\ref{Eqn:det(M_2(t))}) or (\ref{Eqn:det_M(t)_3D}) are positive and the convex hull of the original matrix weighted rational B\'{e}zier curve can be computed by the convex hull of the control points of the converted rational B\'{e}zier curve.
\end{conjecture}

Besides a theoretical proof of the conjecture, investigation of some other simple and effective algorithms for computing the convex hulls of matrix weighted rational B\'{e}zier curves is another interesting future work. As products of B-splines can be finally represented by B-splines, see for example~\citep{ChenRiesenfeldCohen09}, matrix weighted NURBS curves can be converted to conventional NURBS curves in the same way as the conversion of matrix weighted rational B\'{e}zier curves. Alternatively, a matrix weighted NURBS curve can first be decomposed into piecewise matrix weighted rational B\'{e}zier curves and then be converted to a spline of rational B\'{e}zier curves by the method discussed in this paper.

\bibliographystyle{elsarticle-harv}
\bibliography{MatrixWeighted}




%% \bibitem must have one of the following forms:
%%   \bibitem[Jones et al.(1990)]{key}...
%%   \bibitem[Jones et al.(1990)Jones, Baker, and Williams]{key}...
%%   \bibitem[Jones et al., 1990]{key}...
%%   \bibitem[\protect\citepauthoryear{Jones, Baker, and Williams}{Jones
%%       et al.}{1990}]{key}...
%%   \bibitem[\protect\citepauthoryear{Jones et al.}{1990}]{key}...
%%   \bibitem[\protect\astroncite{Jones et al.}{1990}]{key}...
%%   \bibitem[\protect\citepname{Jones et al., }1990]{key}...
%%   \harvarditem[Jones et al.]{Jones, Baker, and Williams}{1990}{key}...
%%


% \end{thebibliography}

\end{document}